\documentclass{article}

\newtheorem{thm}{Theorem}[section]

\newtheorem{prop}[thm]{Proposition}

\newtheorem{lem}[thm]{Lemma}

\newcommand{\be}{\begin{equation}}
\newcommand{\ee}{\end{equation}}
\newcommand{\ben}{\begin{enumerate}}
\newcommand{\een}{\end{enumerate}}
\newcommand{\beq}{\begin{eqnarray}}
\newcommand{\eeq}{\end{eqnarray}}
\newcommand{\beqn}{\begin{eqnarray*}}
\newcommand{\eeqn}{\end{eqnarray*}}
\newcommand{\e}{\varepsilon}
\newcommand{\pa}{\partial}

\newcommand{\qed}{\hspace*{\fill}Q.E.D.}  %Use at end of proof

\title{Einstein  Finsler Metrics and Killing Vector Fields on Riemannian Manifolds  }
\author{Xinyue Cheng\footnote{supported by the National Natural Science Foundation of China (11371386) and the European Union's Seventh Framework Programme (FP7/2007-2013) under grant agreement no. 317721}\\
School of Mathematics and Statistics \\
Chongqing University of Technology \\
Chongqing 400054,  P.R. China \\
Email:  chengxy@cqut.edu.cn\\
\\
 Zhongmin Shen\footnote{supported in part by a NSF grant (DMS-0810159)}\\
 Department of Mathematical Science \\
Indiana University-Purdue University at Indianapolis \\
Indianapolis, USA \\
Email: zshen@math.iupui.edu\\
}

\date{}

\begin{document}
\maketitle

\begin{abstract}

 In this paper, we use a Killing form  on a Riemannian manifold to construct a class of Finsler metrics. We find equations  that characterize Einstein metrics among this class. In particular,  we construct a family of Einstein metrics on $S^3$ with ${\rm Ric} = 2 F^2$,  ${\rm Ric}=0$ and ${\rm Ric}=- 2 F^2$, respectively. This family of metrics provide an important class of Finsler metrics in dimension three, whose Ricci curvature is a constant, but the flag  curvature is not.


\end{abstract}

\section{Introduction}

One of the most important problems in Finsler geometry is to study and characterize Einstein Finsler metrics. By definition, a Finsler metric $F=F(x,y)$ on an $n$-dimensional manifold $M$ is {\it of isotropic Ricci curvature} if
\be
{\rm Ric }= (n-1) K F^2,\label{Einstein}
\ee
where $K= K(x)$ is a scalar function on $M$. Recently, a new notion of Ricci (curvature) tensor  ${\rm Ric}_{ij} $ has been  studied \cite{LiSh} (See Section 2 below for the definition). A Finsler metric $F=F(x,y)$ is {\it of isotropic Ricci tensor}  if
\be
{\rm Ric}_{ij}= (n-1) K g_{ij}, \label{Einstein2}
\ee
where $K=K(x)$ is a scalar function on $M$.

There is an important non-Riemannian quantity $H_{ij}= E_{ij|m}y^m$ defined as the covariant derivative of the mean Berwald curvature $E_{ij}$ along a geodesic. It relates the Ricci curvature tensor $Ric_{ij}$ and the Ricc curvature
$Ric$ as follows
\be
H_{ij}= {\rm Ric}_{ij} - \frac{1}{2} [{\rm Ric}]_{y^iy^j}.  \label{Einstein3}
\ee
Therefore (\ref{Einstein}) and (\ref{Einstein2}) are equivalent for Finsler metrics with $H_{ij}=0$.

A Finsler metric is called an {\it Einstein  metric} if it is of isotropic Ricci curvature and a {\it strong Einstein metric} if it is of {\it isotropic Ricci curvature tensor}. A famous question asked by S. S. Chern is:  whether or not does every smooth manifold admit an Einstein  Finsler metric with Ricci-constant?

Clearly, in dimension two, a Finsler metric is of isotropic Ricci curvature (tensor) if and only if it is of isotropic flag curvature. In dimension three, a Riemannian metric is of isotropic Ricci curvature (tensor) if and only if it is of isotropic sectional curvature, in this case, the sectional curvature must be a constant.   Bao-Robles have proved that the conclusion as above still hold for Randers metrics, that is, a Randers  metric on a manifold of dimension three is of isotropic Ricci curvature (tensor) if and only if it is of constant flag curvature  (\cite{BaRo}).  A natural question arises: is there any Finsler metric with isotropic Ricci curvature (tensor) but not isotropic flag curvature?  The answer is yes. See the following

\begin{thm}\label{thmn=3}  There are Einstein Finsler metrics on $S^3$ with $K=1$, $K=0$ and $K=-1$, respectively.  The metrics take the following form
$F = \alpha\phi(\beta/\alpha)$, where
$\alpha=\sqrt{a_{ij}y^iy^j}$ is a Riemannian metric  and $\beta=b_i(x)y^i$ is a Killing $1$-form of constant length $b$ on $S^3$ satisfying
\[  \overline{\rm Ric}= 2 \alpha^2 -4 (b^2\alpha^2-\beta^2), \ \ \
 s_{0m}s^m_{\ 0} = - (b^2\alpha^2-\beta^2), \ \ \ s^m_{\ 0;m} = 2 \beta   \]
 and $\phi=\phi(s)$ satisfies
\be
(b^2-s^2) [ (1+sQ) Q_s -Q^2+2]-[2sQ+b^2Q^2+1] + K \phi^2 =0,
\ee
where $Q:= \phi'/(\phi-s\phi_s)$,  $\overline{\rm Ric}$ denotes the Ricci curvature of $\alpha$ and $s_{ij} = \frac{1}{2} ( b_{i;j}-b_{j;i})$ the anti-symmetric part of the covariant derivative of $\beta$ with rexpect to $\alpha$.  These metrics are not of constant flag curvature in general, except for the case when $K=1$ and $\phi =1+s$.
\end{thm}

The  Finsler metrics that we constructed on $S^3$  are  called $(\alpha,\beta)$-metrics.
Randers metrics  $F=\alpha+\beta$  are among the simplest non-Riemannian Finsler metrics.
In \cite{BaRo}, Bao-Robles find two equations on $\alpha$ and $\beta$ that characterize Einstein Randers metrics. Using the navigation idea, they can actually classify Einstein Randers metrics upto the classification of Einstein Riemann metrics and the homothetic vector fields.
To search for new Einstein metrics, we shall consider  {\it almost regular} $(\alpha,\beta)$-metrics (see Section 2 below for definition).
Let $\alpha =\sqrt{a_{ij}(x)y^iy^j}$ and $\beta= b_i(x)y^i$.  Put
\[
{r}_{ij}:=\frac{1}{2}({{b}_{i;j}}+{{b}_{j;i}}),\ \ \ {s}_{ij}:=\frac{1}{2}({{b}_{i;j}}-{{b}_{j;i}}),
\]
where ``  ;  " denotes the covariant derivative with respect to the Levi-Civita connection of $\alpha $. Further,  put
\[
{{r}_{j}}:={{b}^{m}}{{r}_{mj}}, \ \ \ \
{{s}_{j}}:={{b}^{m}}{{s}_{mj}},
\]
where $\left( {{a}^{ij}} \right):={{\left( {{a}_{ij}} \right)}^{-1}}$ and $ \ b^{i}:=a^{ij}b_{j}$.   We will denote ${r}_{i0}:={{r}_{ij}}{{y}^{j}},\ {{s}_{i0}}:={{s}_{ij}}{{y}^{j}}$ and
$ {r}_{00}:={{r}_{ij}}{{y}^{i}}{{y}^{j}},\ {{r}_{0}}:={{r}_{i}}{{y}^{i}}, \ {{s}_{0}}:={{s}_{i}}{{y}^{i}}$, etc.  Clearly,  $\beta$ is a Killing form if and only if $r_{ij}=0$.  Thus  $\beta$ is a Killing form of constant length with respect to $\alpha$ if and only if it satisfies the following equations:
\be
r_{ij}  =  0 , \ \ \ s_j  =  0 . \label{EQ2}
\ee
The first author has proved that a regular $(\alpha,\beta)$-metric of non-Randers type is of isotropic S-curvature if and only if (\ref{EQ2}) holds (\cite{Cheng}). Thus (\ref{EQ2})  implies that the S-curvature is a constant, hence $H_{ij}=0$.

 For an $(\alpha,\beta)$-metric $F = \alpha \phi(\beta/\alpha)$ with parallel $\beta$ (i.e. $b_{i;j}=0$), it is an Einstein metric if and only if it is Ricci-flat, and if and only if $\alpha$ is Ricci-flat , regardless of the choice of $\phi$ (see \cite{ChShTi}).

In the following we shall consider  $(\alpha,\beta)$-metrics with non-parallel Killing form $\beta$ of constant length. Here
we allow  the metric to be almost regular in the above sense.
We prove the following

\begin{thm}\label{mainthm}  Let $F =\alpha \phi (s)$, $s=\beta/\alpha$, be an almost regular $(\alpha, \beta)$-metric on an $n$-dimensional manifold. Suppose that $\beta$ is  a non-parallel Killing form with non-zero constant length $b := \|\beta_x\|_{\alpha}$.  Then $F$ is an (strong) Einstein metric if and only if one of the following cases occurs:
\ben
\item[{\rm (i)}] $F$ is an (strong) Einstein $(\alpha, \beta)$-metric of Randers type.
\item[{\rm  (ii)}]
 $\alpha$ and $\beta$ satisfy
\beq
\overline{\rm Ric} & = &  (n-1)\tau \{ K_1 \alpha^2  + K_2 (b^2 \alpha^2-\beta^2)\} - K_3 s_{0m}s^m_{\ 0} ,\label{EQ1*}\\
s^{i}_{\ m}s^m_{\ i} & = & (n-1)\tau,\label{EQ4*}\\
s^m_{\ 0;m} &  = & -(n-1) b^{-2} \tau \beta  \label{EQ5*}
\eeq
and  $\phi=\phi(s)$ satisfies
\be
Q = -b^{-2} s \pm \sqrt{\delta} \sqrt{ 1-b^{-2} s^2 } , \label{ODE*}
\ee
where $Q:=\phi'/(\phi-s\phi')$,  $\delta>0$ is  a constant, $\tau=\tau(x)$ is a scalar function and
\be K_1 = -b^{-2}, \ \  K_2 = b^{-2} ( b^{-2} +\delta), \ \ K_3 =-2 (b^{-2}+\delta).\label{K(1)}
\ee
In this case, ${\rm Ric} =0$ and $F$ is an $(\alpha,\beta)$-metric of Randers type with singularuty.
\item[{\rm (iii)}]  $\alpha$ and $\beta$ satisfy
\beq
\overline{\rm Ric} & = &(n-1) \tau \{  K_1 \alpha^2  + K_2 (b^2 \alpha^2-\beta^2)\},\label{EQ1}\\
s_{0m}s^m_{\; 0} & = & b^{-2} \tau (b^2\alpha^2-\beta^2),\label{EQ4}\\
s^m_{\ 0;m} &  = &  -(n-1) b^{-2} \tau \beta   \label{EQ5}
\eeq
and
$\phi=\phi(s)$ satisfies
\beq
&&  \frac{1}{n-1} (b^2-s^2)[ 2(1+sQ)Q_s-2Q^2]-(2sQ+b^2Q^2)  \nonumber\\
&&+ k [  \delta_1   s^2 + (b^2-s^2) ]  + \delta_2  s^2  + (b^2-s^2)\delta_3  = k b^2 \phi^2,\label{ODE}
\eeq
where   $\delta_1 , \ \delta_2 , \ \delta_3,\  k $  are constants,  $\tau =\tau (x)(\not= 0)$ is a scalar function and
\be
K_1 =  \delta_1 k + \delta_2 , \ \ K_2 = b^{-2} [ (1-\delta_1)k + \delta_3-\delta_2 ] .\label{K(2)}
\ee
In this case, ${\rm Ric} =(n-1)K F^2$ with $K = k\tau $ and $F$ can be of Randers type.
\een
\end{thm}

\bigskip

We are primarily concerned about the existence of $\alpha$ and $\beta$  satisfying the conditions  in Theorem \ref{mainthm}.  At this moment we have not found  any pair $(\alpha,\beta)$  satisfying  (\ref{EQ2}), (\ref{EQ1*}), (\ref{EQ4*}) and (\ref{EQ5*}) with $K_3 \not =0$ and $\tau\not=0$. However,  by using Maple program, we can get a solution of (\ref{ODE*}) as follows:
\[
\phi (s)=C_{1}\phi (s)=C_{1} (\sqrt{b^{2}-s^{2}}\pm \sqrt{\delta} bs ).
\]
In this case, $F=C_{1}( \sqrt{b^{2}\alpha ^{2}- \beta ^{2}}\pm \sqrt{\delta} b\beta )$ is an $(\alpha,\beta)$-metric of Randers type. Obviously, the Riemannian metric $\sqrt{b^{2}\alpha ^{2}- \beta ^{2}}$ is singular at any point $(x,y)\in TM$ satisfying $b\alpha = \beta$.

In Theorem \ref{mainthm} (iii), if we let $c_1 := K_1+K_2b^2$ and $ c_2 := - K_2$ with $\delta_2 =- b^{-2}$ and $k=0$, then (\ref{EQ1}), (\ref{EQ4}) and (\ref{ODE}) are listed in \cite{SeUl} as sufficient conditions for $F$ being  Ricci-flat. As shown in \cite{SeUl}, the  equation (\ref{EQ5}) can be derived from (\ref{EQ1}) and (\ref{EQ4}).

There are many pairs $(\alpha,\beta)$ satisfying (\ref{EQ2}), (\ref{EQ1}), (\ref{EQ4}) and (\ref{EQ5}). We can construct them using Einstein Randers metrics. Let $\bar{F}=\alpha+\beta$ be an Einstein Randers metric with
${\rm Ric}_{\bar{F}} = (n-1) \sigma \bar{F}^2$. If $\beta$  is a Killing form of constant length $b$ satisfying  $s_{0m}s^m_{\ 0} = -\sigma (b^2\alpha^2-\beta^2)$, then, by (7.16) and (7.29) in \cite{ChSh1},  $\alpha$ and $\beta$ satisfy (\ref{EQ2}), (\ref{EQ1}), (\ref{EQ4}) and (\ref{EQ5}) with
\[  \tau = - b^2\sigma, \ \ \ \ K_1 =-b^{-2},\ \ \ \  K_2 = \frac{n+1}{n-1} b^{-2} .\]
Taking constants $\delta_1, \delta_2 , \delta_3, k$ in (\ref{K(2)}) such that  $K_1=-b^{-2}$ and $K_2 =  \frac{n+1}{n-1} b^{-2}$ and letting  $\phi=\phi(s)$ satisfy (\ref{ODE}), we obtain an Einstein metric
$F=\alpha\phi(\beta/\alpha)$ with ${\rm Ric} = (n-1) K F^2$,  where $K= -k \tau$. One can easily verify that $\phi =1+s$ satisfies (\ref{ODE}) with  $ \delta_2 = -b^{-2} (1-\delta_1),$  $ \delta_3 = \frac{n+1}{n-1}$ and $ k= -b^{-2}$. On the other hand, it is surprising that,  for a famous Finsler metric with singularity -- Kropina metric in the form $F=\alpha ^{2}/\beta$, its $\phi (s)=\frac{1}{s}$ satisfies differential equation (\ref{ODE}) with $k=-\frac{1}{4}, \ \delta _{2}=\frac{\delta _{1}}{4}-\frac{1}{b^2}, \   \delta _{3}=\frac{1}{4}-\frac{1}{b^2}.$

The paper is organized as follows. In Section 2,  we give some definitions and notations which are necessary for the present paper, and some lemmas and propositions are contained. The proof of the sufficient condition in Theorem \ref{mainthm} is also given in this section.  Then the proof of the necessary condition in Theorem \ref{mainthm} is given in Section 3.  Further, based on Theorem \ref{mainthm},   we construct Einstein $(\alpha,\beta)$-metrics with $Ric =2F ,  Ric =0$ and $Ric =-2F$ on $S^3$, respectively in Section 4. Here we view $S^3$ as a Lie group. We prove that these constructed metrics are not of constant flag curvature, which means that Theorem \ref{thmn=3} holds. Finally,  we give the conclusion and describe a related reseaech for our research work in Section 5.

\section{Preliminaries}\label{section2}

Let $F=F(x,y)$ be a Finsler metric on a  manifold $M$. The geodesic coefficients $G^{i}$ of $F$ are given  by
\be\label{Gi}
{{G}^{i}}=\frac{1}{4}{{g}^{il}}\left\{ {{\left[ {{F}^{2}} \right]}_{{{x}^{k}}{{y}^{l}}}}{{y}^{k}}-{{\left[ {{F}^{2}} \right]}_{{{x}^{l}}}} \right\},
\ee
where $g_{ij}(x,y):=\frac{1}{2}\big[F^{2}\big]_{y^{i}y^{j}}(x,y)$ and $(g^{ij}):=(g_{ij})^{-1}$.

The well-known non-Riemannian quantity,   S-curvature,  is given by
\be
{\bf S} = \frac{\partial G^m}{\partial y^m}- y^m \frac{\partial }{\partial x^m} [\ln \sigma_F ],
\ee
where $dV_F =\sigma_F(x) dx^1\cdots dx^n$ is the Busemann-Hausdorff volume form.
By definition,  $F$ is said to be of {\it isotropic S-curvature} if there exists  a scalar function $c(x)$ on $M$ such that ${\bf S}(x,y)=(n+1) c(x)F(x,y)$. If $c(x)=$constant, we say that $F$ has   {\it constant S-curvature}.

Using the S-curvature, one obtains another quantity $H= H_{ij} dx^i\otimes dx^j$  given by
\[  H_{ij} =\frac{1}{2} {\bf S}_{ij; m }y^m.\]

For any $x\in M$ and $y\in T_{x}M\backslash\{0\}$, the {\it Riemann curvature} ${\bf R}_{y} = R^i_{\ k} \frac{\pa }{\pa x^i} \otimes dx^k$ of $F$ is defined by
\be
R^i_{\ k} = 2 \frac{\pa G^i}{\pa x^k} -\frac{\pa^2 G^i}{\pa x^m\pa y^k}y^m + 2 G^m \frac{\pa^2 G^i}{\pa y^m \pa y^k} -\frac{\pa G^i}{\pa y^m} \frac{\pa G^m}{\pa y^k}. \label{Rik}
\ee
The Riemann curvature tensor $R^{\ i}_{j \ kl}$ is given by
\[
R^{\ i}_{j \ kl}=\frac{1}{3}\left\{R^{i}_{\ k \cdot l \cdot j}-R^{i}_{\ l \cdot k \cdot j} \right\}.
\]
It is easy to see that $R^{i}_{\ k}=R^{\ i}_{j \ kl}y^{j}y^{l}$. The two curvature tensors $R^{\ i}_{j \ kl}$ and $R^{i}_{\ k}$ essentially contain the same geometric data.
The {\it Ricci curvature tensor} ${\rm Ric}_{ij}$ is defined by
\be
{\rm Ric}_{ij}:=\frac{1}{2}\left\{R^{\ m}_{i \ m j}+R^{\ m}_{j \ m i} \right\}.
\ee
The {\it Ricci curvature} is the trace of the Riemann curvature, which is defined by
\be
{\rm Ric}=R^{m}_{\ m}.
\ee
We have ${\rm Ric}={\rm Ric}_{ij}y^{i}y^{j}=R^{\ m}_{i \ m j}y^{i}y^{j}$.
 A Finsler metric $F$ on an $n$-dimensional manifold $M$ is called an {\it  Einstein metric} if the Ricci curvature satisfies
\be
{\rm Ric}=(n-1)K F^2, \label{Ric1}
\ee
where $K=K(x)$ is a scalar function on $M$.
$F$ is said to be of {\it Ricci constant} if $F$ satisfies (\ref{Ric1}) and $K =$ constant. We call $F$  the {\it Ricci flat} Finsler metric if ${\rm Ric}=0$.

We  have the following important identity (\cite{LiSh}):
\be
{\rm Ric}_{ij}-\frac{1}{2} [{\rm Ric}]_{y^iy^j}  = H_{ij}.
\ee

\bigskip
We now focus on $(\alpha, \beta)$-metrics.
By  definition, an $(\alpha,\beta)$-metric  on a manifold $M$ is
expressed in the following form
\[
F=\alpha\phi(s),\ s=\frac{\beta}{\alpha},
\]
where $\alpha=\sqrt{a_{ij}(x)y^{i}y^{j}}$ is a Riemannian metric and $\beta= b_{i}(x)y^{i}$ is a 1-form on $M$. It is proved (\cite{BaChSh}) that  $F = \alpha \phi(\beta/\alpha)$ is a positive definite Finsler
metric for any $\alpha, \beta$ with $\|\beta\|_{\alpha}(x) < b_o, \ x\in M$, if and only if the function $\phi = \phi(s)$ is a $C^{\infty}$ positive function on an open interval $(-b_{0}, b_{0})$ satisfying
\[
\phi(s)-s\phi'(s)+(b^{2}-s^{2})\phi''(s)>0,\ \ |s|\leq b<b_{0}.
\]
For the above function $\phi =\phi (s)$, if the 1-form $\beta$ satisfies $b(x):=\|\beta\|_{\alpha}(x) \leq b_{0}$, then $F=\alpha \phi (\beta/\alpha)$ might be singular at the points $x\in M$ with $b(x)=b_{0}.$  Such metrics are called {\it almost regular $(\alpha, \beta)$-metrics}.
When $\phi =1+s$, the Finsler metric $F=\alpha+\beta$ is just the  Randers metric. When $\phi=1/ s$, the Finsler metric $F={\alpha^2}/{\beta}$ is called {\it Kropina metric}.  Randers metrics and Kropina metrics are both C-reducible (\cite{Ma1}). However, Randers metrics are regular Finsler metrics but Kropina metrics are Finsler metrics with singularity.

We have the following

\begin{lem}\label{aGi}{\rm (\cite{BaChSh}\cite{ChernSh})} For an $(\alpha ,\beta )$-metric $F=\alpha \phi (s),s=\beta /\alpha $,  the geodesic coefficients ${{G}^{i}}$ of $F$ and the geodesic coefficients $\bar{G}^i$ of $\alpha$ are related by
\be
{{G}^{i}}=\bar{G}+\alpha Q{{s}^{i}}_{0}+\alpha^{-1}\Theta \left\{ {{r}_{00}}-2Q\alpha {{s}_{0}} \right\}y^i+\Psi \left\{ {{r}_{00}}-2Q\alpha {{s}_{0}} \right\}{{b}^{i}},  \label{Giofab}
\ee
where
\beqn
 Q&:=&\frac{{{\phi }'}}{\phi -s{\phi }'}, \\
 \Theta &:=& \frac{\phi {\phi }'-s(\phi {\phi }''+{\phi }'{\phi }')}{2\phi \left[ (\phi -s{\phi }')+(b^2-{{s}^{2}}){\phi }'' \right]}, \\
\Psi &:=&\frac{{{\phi }''}}{2\left[ (\phi -s{\phi }')+(b^2-{{s}^{2}}){\phi }'' \right]}.
\eeqn
\end{lem}

In \cite{ChSh0}, we have proved that a Randers metric $F=\alpha +\beta$ is of isotropic S-curvature, ${\bf S}=(n+1)c(x)F$, if and only if
\be
r_{ij}+b_{i}s_{j}+b_{j}s_{i}=2c(x)(a_{ij}-b_{i}b_{j}). \label{RandersS}
\ee
Further, we considered Finsler metrics of Randers type in the following form
$$
F=k_{1}\sqrt{\alpha ^2 +k_{2}\beta ^2}+k_{3}\beta ,
$$
where $k_{1}>0, k_{2}$ and $k_{3}\neq 0$ are constants. We obtained the sufficient and necessary condition that a Finsler metric of Randers type to be of isotropic S-curvature (\cite{ChSh}). More general, we characterized almost regular $(\alpha, \beta)$-metrics of non-Randers type with isotropic S-curvature (\cite{ChSh}). For regular $(\alpha,\beta)$-metrics, the first author has proved the following theorem.

\begin{thm} {\rm (\cite{Cheng})} \label{mainthonS}
A regular $(\alpha,\beta)$-metric $F=\alpha \phi (\beta/\alpha)$ of non-Randers type on an $n$-dimensional  manifold $M$ is of isotropic S-curvature, ${\bf S}=(n+1)cF$, if and only if $\beta$ satisfies
\be
r_{ij}=0, \ \ \ s_{j}=0. \label{simplification}
\ee
In this case, ${\bf S}=0$, regardless of the choice of a particular $\phi =\phi (s)$.
\end{thm}

For an almost regular $(\alpha,\beta)$-metric, if (\ref{simplification}) holds, then ${\bf S}=0$. The converse might not be true.

We need the following lemma.

\begin{lem}  Assume that (\ref{simplification}) holds, then
\be
b^m s_{jm;k} = - s_{jm} s^m_{\ k}.
\ee
In particular, we have
\be
 b^i s^m_{\ i;m} = - s^i_{\ m}s^m_{\ i}. \label{sss}
\ee
\end{lem}
{\it Proof}:
\begin{eqnarray*}
0 & = &  s_{j;k} = b^m_{\ ; k} s_{mj} + b^m s_{mj;k}\\
& = & b_{m;k} s^m_{\ j} - b^m  s_{jm;k}\\
& = & s_{mk}s^m_{\ j} -b^m s_{jm;k}\\
& = & - s_{jm} s^m_{\ k} - b^m s_{jm;k}.
\end{eqnarray*}
\qed

\bigskip

Let $F = \alpha \phi \Big (\frac{\beta}{\alpha}\Big )$ be an $(\alpha ,\beta )$-metric. In the following,
we always assume that $\beta $ is a Killing form of constant length $b$, namely it satisfies (\ref{EQ2}).
\be
r_{ij}  =  0 , \ \ \ s_j  =  0 . \label{EQ2**}
\ee
By (\ref{Giofab}), we can rewrite the geodesic coefficients of an $(\alpha,\beta)$-metric as
\be
{{G}^{i}}=\bar{G}+{{T}^{i}}, \label{GGab}
\ee
where ${{T}^{i}}:=\alpha Q{{s}^{i}}_{0}$.  From (\ref{Rik}) and by use of a technique for computing Riemannian curvature, we have
\[
{{R}^{i}}_{j}=\bar{R}^i_{\ j}+R{{T}^{i}}_{j},
\]
where $\bar{R}^i_{\ j}$ denote the Riemann curvature of $\alpha$ and
\[
R{{T}^{i}}_{j}:=2{{T}^{i}}_{;j}-{{({{T}^{i}}_{;k})}_{{{y}^{j}}}}{{y}^{k}}+2{{T}^{k}}{{({{T}^{i}})}_{{{y}^{j}}
{{y}^{k}}}}-{{({{T}^{i}})}_{{{y}^{k}}}}{{({{T}^{k}})}_{{{y}^{j}}}}.
\]
After a series of complex computations and  by use of Maple program, we can obtain

\begin{prop}{\rm  (\cite{ChShTi})}\label{RiemCur}
For any $(\alpha ,\beta )$-metric $F=\alpha \phi (s),\ s=\beta /\alpha $, the Riemannian curvature is given by
\be
{{R}^{i}}_{j}=\bar{R}^i_{\ j}+R{{T}^{i}}_{j},  \label{abRiemannCur}
\ee
where
\beqn
{RT}^{i}_{\ j}&=& C_{331}s^{i}_{\ 0}s_{0j} + \alpha C_{332}s^{i}_{\ k}s^{k}_{\ 0}l_{j}+\alpha C_{333}s^{i}_{\ k}s^{k}_{\ 0}b_{j}-(Q^{2}\alpha ^2)s^{i}_{\ k}s^{k}_{\ j}\\
&&+(2Q\alpha)s^{i}_{\ 0|j}-(Q\alpha)s^{i}_{\ j|0}+[C_{311}s^{i}_{\ 0|0}l_{j}-Q_{s}s^{i}_{\ 0|0}b_{j}]
\eeqn
and
\beqn
C_{331}&=&-3Q^{2}+3sQQ_{s}+3Q_{s},\\
C_{340}&=&-2[Q_{s}(B-s^2)+sQ+1]\Psi+Q_{s},\\
C_{332}&=&(Q-sQ_{s})Q, \\
C_{333}&=& QQ_{s},\\
C_{311}&=& sQ_{s}-Q.
\eeqn
\end{prop}

By Proposition 3.3 in \cite{ChShTi}, we  get the following

\begin{prop}\label{proRic} For an $(\alpha ,\beta )$-metric $F=\alpha \phi (s),s=\beta /\alpha $, if (\ref{simplification}) holds, then the Ricci curvature   of $F$ is related to the Ricci curvature $\overline{\rm Ric}$ of $\alpha $ by
\be
{\rm Ric}= \overline{\rm Ric}+ s_{0m}s^m_{\ 0} c_{19} +\alpha s^m_{\ 0;m} c_{24} +\alpha^2 s^i_{\ m}s^m_{\ i} c_{26}, \label{eqRic}
\ee
where
\beqn
c_{19}&=& -2Q^{2}+2(1+sQ)Q_{s},\\
c_{24}&=& 2Q, \\
c_{26}&=& -Q^{2}.
\eeqn
\end{prop}

\bigskip
\noindent
{\it Proof of the sufficient condition in Theorem \ref{mainthm}}:

(ii)  Assume that $\alpha$ and $\beta$ satisfy (\ref{EQ1*})-(\ref{EQ5*}) with
\[ K_1 =  -b^{-2}, \ \ \ \ \ K_2 = b^{-2}(b^{-2}+\delta), \ \ \ \ K_3 =-2(b^{-2}+\delta) \]
and $\phi$ satisfies (\ref{ODE*}).
Then  $Q $ satisfies the following two equations:
\beqn
&& c_{19} = K_3 ,\\
&& K_1 + K_2 (b^2-s^2)- 2 b^{-2} s Q - Q^2 =0.
\eeqn
Then
it follows from (\ref{eqRic}) that
\begin{eqnarray*}
 {\rm Ric  }&  = &  (n-1)\tau \{ K_1 \alpha^2+K_2 (b^2\alpha^2 -\beta^2) \}- K_3 s_{0m}s^m_{\ 0} + c_{19}s_{0m}s^m_{\ 0}  \\
&& -(n-1)b^{-2} \tau  c_{24}\alpha\beta + (n-1)\tau c_{26} \alpha^2  \\
&= & (n-1) \tau \alpha^2 \Big \{  K_1 + K_2 (b^2-s^2) - b^{-2}  s  c_{24} + c_{26}  \Big \}
=0.
\end{eqnarray*}
This proves the sufficient condition in Theorem \ref{mainthm} (ii).

(iii)  Assume that $\alpha$ and $\beta$ satisfy (\ref{EQ1})-(\ref{EQ5})
with $K_1, K_2$ given by
\[ K_1 = \delta_1 k +\delta_2 , \ \ \ \ \
K_2 = b^{-2} [(1-\delta_1)k + (\delta_3-\delta_2)] ,\]
and $\phi$ satisfies (\ref{ODE}).
Then it follows from (\ref{eqRic}) that
\begin{eqnarray*}
 {\rm Ric } &  = & (n-1)\tau \{  K_1 \alpha^2+K_2 (b^2\alpha^2 -\beta^2)\} \\
&& + c_{19} b^{-2} \tau (b^2\alpha^2 -\beta^2)  -(n-1) c_{24} b^{-2} \tau  \alpha\beta + (n-1) c_{26}  \tau \alpha^2\\
& = &\Big \{ (n-1)b^2[ k+\delta_3  ] - (n-1)[ (1-\delta_1) k +(\delta_3-\delta_2)] s^2 \\
&&  + c_{19}   (b^2 -s^2) -(n-1) c_{24}    s  + (n-1) c_{26} b^2 \Big \}b^{-2}\tau \alpha^2\\
&= &(n-1)\Big \{  k [ \delta_1 s^2 + (b^2-s^2)]   + \delta_2 s^2 +(b^2-s^2) \delta_3   \\
&& +\frac{1}{n-1} (b^2-s^2) c_{19}-( c_{24}s- b^2 c_{26} ) \Big \} b^{-2} \tau \alpha^2\\
& =& (n-1) k  \tau \phi^2 \alpha^2 =(n-1)  k\tau  F^2.
\end{eqnarray*}
This proves the sufficient condition in Theorem \ref{mainthm} (iii).

\section{Proof of Theorem \ref{mainthm}}

In this section, we shall prove the necessary condition in Theorem \ref{mainthm}.

Let $F=\alpha \phi(s), \ s= \beta/\alpha$, be an $(\alpha ,\beta )$-metric on an $n$-dimensional  manifold $M$. Firstly, we take an orthonormal basis on $T_{x}M$  with respect to $\alpha $ at any fixed point $x$ such that
\be
\alpha =\sqrt{\sum{{{({{y}^{i}})}^{2}}}},\ \ \ \beta =b{{y}^{1}}. \label{specialframe}
\ee
To simplify the computations, we   take the following coordinate transformation
$\psi:  (s,u^A) \to (y^i)$:
\be
y^1 = \frac{s}{\sqrt{b^2-s^2} }\bar{\alpha}, \ \ \ \ \ \ y^A = u^A,\label{specialcoordinates}
\ee
where $\bar{\alpha} =\sqrt{\sum_{A=2}^n (u^A)^2}$.  Here, our index conventions are  as follows:
\[
1\leq i,j,k, \cdots \leq n, \ \ \ \ 2\leq A, B, C, \cdots \leq n.
\]
Then
\be
\alpha =\frac{b}{\sqrt{b^2-s^2} }\bar{\alpha}, \ \ \ \ \
\beta = \frac{bs}{\sqrt{b^2-s^2} } \bar{\alpha}. \label{abara}
\ee
Thus
\[
F = \alpha \phi (\beta/\alpha) = \frac{b\phi(s)}{\sqrt{b^2-s^2}} \bar{\alpha}.
\]
Further,
\beqn
s_{0m}s^{m}_{\ 0}&=&\frac{{{s}^{2}}{{s}_{1m}}s^{m}_{\ 1}}{{{b}^{2}}-{{s}^{2}}}{{{\bar{\alpha }}}^{2}}+\frac{2{{s}_{1m}}\bar{s}^{m}_{\ 0}}{\sqrt{{{b}^{2}}-{{s}^{2}}}}s\bar{\alpha }+{{{\bar{s}}}_{0m}}\bar{s}^{m}_{\ 0}, \\
s^{m}_{\ 0;m}&=&\frac{s^{m}_{\ 1;m}}{\sqrt{{{b}^{2}}-{{s}^{2}}}}s\bar{\alpha }+\bar{s}^{m}_{\ 0;m},
\eeqn
where $\bar{s}^{m}_{\ 0}:=s^{m}_{ \ A}y^{A}, \ \bar{s}_{0m}\bar{s}^{m}_{\ 0}:=s_{Am}s^{m}_{ \ B}y^{A}y^{B}$, etc..

Since $\bar{G}=\frac{1}{2}{}^{\alpha }\Gamma _{jk}^{i}(x){{y}^{j}}{{y}^{k}}$, where ${}^{\alpha }\Gamma _{jk}^{i}$ denote the Christoffel symbols  of  $\alpha $, the Ricci curvature $\overline{\rm Ric} = G_{ij} (x)y^iy^j$ is quadratic polynomial in $y$. Hence, it is easy to prove the following lemma.

\begin{lem}\label{lem4.2}
The Ricci curvature of $\alpha$ can be expressed as
\[
\overline{\rm Ric}=\frac{{{s}^{2}}{{G}_{11}}}{{{b}^{2}}-{{s}^{2}}}{{\bar{\alpha }}^{2}}+({{\bar{G}}_{10}}+{{\bar{G}}_{01}})\frac{s}{\sqrt{{{b}^{2}}-{{s}^{2}}}}\bar{\alpha }+{{\bar{G}}_{00}}.
\]
\end{lem}

Assume that $F=\alpha\phi(\beta/\alpha)$ is an Einstein metric, ${\rm Ric}=(n-1)K F^2$, and satisfies (\ref{EQ2}). Then,  by (\ref{eqRic}) and Lemma \ref{lem4.2}, we obtain
\be
\Xi _{4}{{\bar{\alpha }}^{2}}+{{\Xi }_{3}}{{\bar{\alpha }}}+{{\Xi }_{2}}=0, \label{keyeq}
\ee
where
\beqn
{{\Xi }_{4}}&=&{{G}_{11}} s^2+s^{m}_{\ 1;m}{{c}_{24}}bs+{b^2} s^{i}_{\ m}s^{m}_{\ i}{{c}_{26}} -(n-1)K {{b}^2}{{\phi }^{2}}, \\
 {{\Xi }_{3}}&=&\sqrt{{{b}^{2}}-{{s}^{2}}}\left\{
({{{\bar{G}}}_{10}}+{{{\bar{G}}}_{01}}) s+\bar{s}^{m}_{\ 0;m}{{c}_{24}}b \right\}, \\
{{\Xi }_{2}}&=&({{b}^{2}}-{{s}^{2}})\left\{{{{\bar{s}}}_{0m}}\bar{s}^{m}_{\ 0}{{c}_{19}}+{{{\bar{G}}}_{00}} \right\} . \\
\eeqn
From (\ref{keyeq}),
we obtain the following fundamental equations
\beq
&& \Xi _{4}{{\bar{\alpha }}^{2}}+{{\Xi }_{2}}=0 , \label{X_4} \\
&& \Xi_3=0.  \label{X_3}
\eeq

We may assume that $c_{24}/s=2Q/s \not= constant$.  Otherwise  $\phi(s)= k_1 \sqrt{1+k_2 s^2}$ for some constants $k_1, k_2$, and hence $F$ is Riemannian, which is excluded in the theorem. Then  it follows from (\ref{X_3}) that
\be
  \bar{G}_{10} +\bar{G}_{01} =0, \ \ \ \ \ \ \bar{s}^m_{\ 0;m}=0.  \label{sss2}
\ee

For simplicity, we let
\[  C(s):= \frac{\Xi_4}{b^2-s^2}.\]
Then we can rewrite  (\ref{X_4}) as
\be
 c_{19} \bar{s}_{0m}\bar{s}^m_{\ 0}  +\bar{G}_{00} + C(s) \bar{\alpha}^2 =0  \label{C1}
\ee
 and get
\be
c_{19}' \bar{s}_{0m}\bar{s}^m_{\ 0}  + C'(s) \bar{\alpha}^2=0.  \label{C2}
\ee

\bigskip
\noindent{\bf Case 1}:  $c_{19}'(s)=0$.

By (\ref{C2}), one can easily see that   $c_{19}(s)=C_3$  and $C(s)= C_4$ are independent of $s$. $C_3$ is a constant since $c_{19}$ is independent of $x\in M$ while $C_4$ might be a function of $x\in M$.
Let
\[ C_1 := G_{11}, \ \ \ \  M_1:=s^m_{\ 1;m}, \ \ \ \ M_2:=s^i_{\ m}s^m_{\ i}.\]
By the definitions of $c_{19}(s)$ and $C(s)$, we obtain two ODEs on $\phi$:
\begin{eqnarray}
&& -2 Q^2 + 2 (1+sQ)Q_s = C_3 ,  \label{C_3}\\
&& C_1 s^2 + 2b Q M_1 s - b^2 M_2 Q^2 -(n-1) K  b^2 \phi^2 = C_4 (b^2-s^2).  \label{K_4}
\end{eqnarray}
We rewrite (\ref{C1}) as follows
\be
\bar{G}_{00} =-  C_4 \bar{\alpha}^2 - C_3  \bar{s}_{0m}\bar{s}^m_{\ 0}.
\ee
In general coordinates, we have
\be
\overline{\rm Ric}= C_1 \alpha^2+ C_2 (b^2\alpha^2-\beta^2) - C_3 s_{0m}s^m_{\ 0},
\ee
where $C_2$ is determined by
\be
C_1 + b^2 C_2 =  - C_4 . \label{C22}
\ee
By (\ref{sss2}),
\be
s^m_{\ 0;m}  = b^{-1} M_1 \beta .
\ee
By (\ref{sss}), we have
\be
M_1 = - b^{-1} M_2.    \label{M1M2}
\ee

\bigskip
\noindent{\bf Case 1a}:
If $C_3=0$, then by (\ref{C_3}), we get
\be
Q = \frac{ s\pm \sqrt{k + s^2}}{k},
\ee
where $k>0$ is a constant.
Since
\[ \phi = e^{\int \frac{Q(s)}{1+sQ(s)} ds } .\]
we get
\[ \phi = \frac{\sqrt{k+s^2}\pm s}{\sqrt{k}}.\]
Thus $F=\alpha \phi (\beta/\alpha)$ is an Einstein $(\alpha,\beta)$-metric of Randers type.

\bigskip

\noindent{\bf Case 1b}:
If $C_3 \not=0$, then by (\ref{C_3}), we get
\be
s-\frac{2Q}{C_3} - k \sqrt{2Q^2+C_3}=0,
\ee
where $k$ is a non-zero constant, that is,
\[ Q = \frac{ \frac{2s}{C_3} \pm \sqrt{ \Big ( \frac{2s}{C_3} \Big )^2 - \Big ( \frac{4}{C_3^2} - 2 k^2  \Big ) ( s^2-C_3 k^2)   }  }{ \frac{4}{C_3^2}- 2 k^2}.\]
We can simplify express $Q$ in the following form
\be  Q = \frac{\lambda}{\delta}  s \pm \sqrt{ \delta +  \lambda s^2},\label{QQQ}
\ee
where $\lambda$  and $\delta $  are constants and $\delta >0$. In this case, $C_3= \frac{ 2 (\lambda-\delta^2)}{\delta} $.

If $K \neq 0$, then by (\ref{K_4}), we get
\[  \phi^2 = \frac{ C_1 s^2 + 2b M_1 s Q-b^2 M_2 Q^2 -C_4 (b^2-s^2) }{ (n-1) b^2 K }.\]
Plugging (\ref{QQQ}) into the above expression, one get
\[
\phi^2 = a  + bs^2 + c s \sqrt{ \delta + \lambda s^2},
\]
where $a$, $b$ and $c$ are constants. Using the above formula, we get
\begin{eqnarray*}
 Q  & =  & \frac{ (\phi^2)'}{2\phi^2-s(\phi^2)'}\\
& = & \frac{ c\delta + 2c \lambda s^2 + 2bs \sqrt{\delta + \lambda s^2 }  }{ c\delta s + 2a \sqrt{\delta +\lambda s^2}}.
\end{eqnarray*}
Comparing it with (\ref{QQQ}),
we get
\[  b= a\frac{\lambda +\delta^2}{\delta}, \ \ \ \ c = 2a.\]
Then
\[ \phi^2= a \Big \{ 1 + (\frac{\lambda}{\delta} +\delta) s^2 + 2 s \sqrt{ \delta + \lambda s^2}\Big \} = a \Big ( \sqrt{\delta} s +  \sqrt{ 1+ \frac{\lambda}{\delta} s^2}   \Big )^2.\]
Thus $F$ is of Randers type.

\bigskip

If $K=0$, we are going to  express $M_1, M_2$ and $C_3$ in terms of some constants and $C_4$.  Now (\ref{K_4}) is reduced to
\be
C_1 s^2 + 2b Q M_1 s - b^2 M_2 Q^2 = C_4 (b^2-s^2). \label{C1M1}
\ee
Plugging (\ref{QQQ}) into (\ref{C1M1}) yields
\begin{eqnarray}
&& C_1+2bM_1\frac{\lambda}{\delta}-b^2M_2\lambda-b^2M_2\Big (\frac{\lambda}{\delta}\Big )^2+C_4=0 ,\\
&& M_1-b M_2\frac{\lambda}{\delta}=0 ,\\
&& C_4+M_2\delta=0.
\end{eqnarray}
We get
\[ C_1 = \frac{ b^2 \lambda (\lambda -\delta^2 ) -\delta^3}{\delta^3} C_4, \ \ \ \
M_1 =-\frac{b\lambda}{\delta^2} C_4, \ \ \ \ \ M_2 = -\frac{1}{\delta} C_4.\]
In this case,
\begin{eqnarray}
\overline{\rm Ric} & = & C_1 \alpha^2 + C_2 (b^2\alpha^2-\beta^2) - C_3 s_{0m}s^m_{\ 0}, \label{two1} \\
s^i_{\ m}s^m_{\ i} & = & M_2 =-\frac{1}{\delta}C_{4}, \label{two2}\\
s^m_{\ 0;m} & = &  b^{-1}M_1 \beta =-\frac{\lambda}{\delta ^2}\beta C_{4}, \label{two3}
\end{eqnarray}
where, $C_3 = \frac{2 (\lambda-\delta^2)}{\delta}$ and
\[ C_2 = -\frac{\lambda (\lambda - \delta^2)}{\delta^3}C_4 \]
 by (\ref{C22}). From (\ref{M1M2}), we have
\be
(\delta +b^2 \lambda) C_4 =0.  \label{aK_4}
\ee
By the assumption, $s_{ij}\not=0$. Thus the matrix $(s^i_{\ m}s^m_{\ j})$ is semi-negative definite. Then $M_2 = s^{i}_{\ m}s^{m}_{\ i} <0$, which means that $C_{4}>0$. It follows from (\ref{aK_4}) that
\[  \lambda = - b^{-2}\delta.\]
Further, we may define a scalar function $\tau =\tau(x)$ such that
\[
C_4 = - (n-1)\delta \tau .
\]
We get
\begin{eqnarray*}
C_1 & = & -(n-1) \tau \frac{ b^2 \lambda (\lambda -\delta^2 ) -\delta^3}{\delta^2}  = - (n-1) \tau b^{-2}, \\
C_2 & = & (n-1)\tau  \frac{\lambda (\lambda - \delta^2)}{\delta^2} = (n-1) b^{-2} ( b^{-2}+\delta)\tau ,\\
C_3 & = & -2 (b^{-2}+\delta).
\end{eqnarray*}
Plugging them into  (\ref{QQQ}), (\ref{two1}), (\ref{two2}) and (\ref{two3}) and letting $K_{1}:=-b^{-2}$, $K_{2}:=b^{-2}(b^{-2}+\delta)$ and $K_{3}:=C_{3}$ yield  Theorem \ref{mainthm} (ii).

\bigskip
\noindent
{\bf Case 2}:
 $c_{19}'(s)\not= 0$ for some $s$.  Then by (\ref{C2}), there is a number $ \tau  =\tau (x)$  independent of $s$ such that
\be
C'(s) = - \tau  c_{19}' . \label{P1}
\ee
Equivalently,  there is a scalar function $\rho=\rho (x)$ independent of $s$ such that
\be
C(s)+\tau  c_{19} =- \rho.  \label{CK}
\ee
By (\ref{C2}) and (\ref{P1}), we obtain
\be
 \bar{s}_{0m}\bar{s}^m_{\ 0} = \tau  \bar{\alpha}^2 .  \label{barss}
\ee
Then (\ref{C1}) becomes
\be
\bar{G}_{00} = \rho  \bar{\alpha}^2. \label{barG00}
\ee
Let us rewrite (\ref{CK}) as
\be
G_{11}s^2 + s^m_{\ 1;m} c_{24} bs +b^{2} s^i_{\ m}s^m_{\ i} c_{26}  - (n-1)K b^2 \phi^2 + (b^2-s^2)(\tau  c_{19} + \rho )=0.  \label{ode2}
\ee

In a general coordinate system,  by Lemma \ref{lem4.2} and (\ref{sss2}), (\ref{barG00}), we can get
\be
\overline{\rm Ric} = C_1 \alpha^2 + C_2 ( b^2 \alpha^2 -\beta^2 ), \label{equaRica}
\ee
where
\be
C_{1}=G_{11}, \ C_{2}=b^{-2}(\rho -G_{11}).\label{C_1C_2}
\ee
  By
(\ref{barss}) and $s_{1m} = b^{-1} s_m =0$, we have the following
\be
s_{0m}s^m_{\ 0}  =  b^{-2} \tau  (b^2 \alpha^2-\beta^2). \label{s0mm}
\ee
It follows from (\ref{s0mm}) that
\be
s^i_{\ m} s^m_{\ i} = (n-1) \tau . \label{ismmi}
\ee
By (\ref{sss}), we get
\be
  s^m_{\ 1;m}  = -b^{-1} s^i_{\ m} s^m_{\ i} =-  (n-1) b^{-1} \tau . \label{sm1m}
\ee
Further, by (\ref{sss2}) and (\ref{abara}), we have
\[
s^{m}_{\ 0;m}=\frac{s^{m}_{\ 1;m}}{\sqrt{b^2 -s^2}}s\bar{\alpha}+\bar{s}^{m}_{\ 0;m}=s^{m}_{\ 1;m}b^{-1}\beta = -(n-1) b^{-2} \tau  \beta.
\]
Then we have the following
\begin{lem} \label{basiclemma1} Assume that an almost regular $(\alpha,\beta)$-metric  $F=\alpha\phi(\beta/\alpha)$ is an Einstein metric, ${\rm Ric}=(n-1)K F^2$, and satisfies (\ref{EQ2}) and $c_{19}' (s)\neq 0$. Then
\begin{eqnarray}
\overline{\rm Ric} & = & C_1 \alpha^2 + C_2 ( b^2 \alpha^2 -\beta^2 ), \label{M1}\\
s_{0m}s^m_{\ 0} & = & b^{-2} \tau  (b^2 \alpha^2-\beta^2),  \label{M2} \\
s^m_{\ 0;m} & = &  - (n-1) b^{-2} \tau  \beta  , \label{M3}
\end{eqnarray}
where $C_{1}=G_{11}, \ C_{2}=b^{-2}(\rho -G_{11})$ and $ \tau  =\tau (x), \ \rho=\rho (x)$ are scalar functions  independent of $s$ on the manifold.
\end{lem}

\bigskip
Now, by (\ref{ismmi}) and (\ref{sm1m}), (\ref{ode2}) becomes
\be
C_1 s^2 - (n-1) \{ \tau  c_{24} s-  \tau  c_{26} b^2 +K b^2 \phi^2\} + (b^2-s^2)(\tau  c_{19} + \rho)=0.  \label{ode2*****}
\ee

\begin{lem} \label{basiclemma2} If (\ref{ode2*****}) holds, then
\be
C_1 = (n-1) (\delta_1 K + \delta_2 \tau) , \label{G11}
\ee
\be
\rho = (n-1)(K + \delta_3  \tau) ,\label{K_2}
\ee
where $\delta_1 , \delta _2 , \delta _3$ are constants.
\end{lem}
{\it Proof}: Firstly, letting $s=0$ in (\ref{ode2*****}) yields  (\ref{K_2}), where
$$
\delta _3 =-\left[c_{26}(0)+ \frac{1}{n-1} c_{19}(0) \right]
$$
is a constant.

Differentiate (\ref{ode2*****}) with respect to $s$ twice. Then, by  (\ref{K_2}) and  letting $s=0$ yields (\ref{G11}) .
\qed

\bigskip

Now (\ref{ode2*****}) become
\beqn
&& \{  \delta_1 K + \delta_2 \tau \} s^2 -  \{ \tau  c_{24} s-  \tau  c_{26} b^2 +K b^2 \phi^2\}\\
&& + (b^2-s^2)\{ \frac{1}{n-1}\tau  c_{19} +K + \delta_3 \tau  \}=0,
\eeqn
that is,
\beq
&&  K \{  \delta_1   s^2- b^2 \phi^2 + (b^2-s^2) \} \nonumber \\
&&  + \tau  \left\{ \delta_2  s^2 +(b^2-s^2) \delta_3 -( c_{24} s -c_{26} b^2) + \frac{1}{n-1}(b^2-s^2)c_{19} \right\}  =0. \label{basiceq23}
\eeq
Note that $K =K (x)$ and $\tau  = \tau (x)$ are scalar functions  on $M$ while
their coefficients are functions of $s$, independent of $x\in M$.

Under the condition that $\beta$ is not parallel with respect to $\alpha$, we assert that $\tau \neq 0$. In fact, if $\tau =0$ at some points, then we can conclude that $\beta$ is closed, $s_{ij}=0$, by (\ref{M2}). Further,  we can see that $\beta$ is parallel with respect to $\alpha$ by (\ref{EQ2}).  It is a contradiction. Thus one can see that $K = k \tau $ for some constant $k$ from (\ref{basiceq23}). Then (\ref{basiceq23}) is reduced to

\beq
&&  k\{  \delta_1   s^2- b^2 \phi^2 + (b^2-s^2) \} \nonumber \\
&&  + \left\{ \delta_2  s^2+ (b^2-s^{2})\delta_3 -( c_{24} s -c_{26} b^2) + \frac{1}{n-1}(b^2-s^2)c_{19}  \right\}  =0. \label{basiceq3}
\eeq
It is just (\ref{ODE}). In this case, by Lemma \ref{basiclemma1} and Lemma \ref{basiclemma2},
\[
C_1 = (n-1)[\delta_1 k+ \delta _2  ] \tau , \ \ \  \
C_2 = (n-1) b^{-2}[(1-\delta_1)k+(\delta_3-\delta_2) ] \tau.
\]
This completes the proof of Theorem \ref{mainthm}  (iii) by letting $K_{1}:=\delta _{1}k+\delta _{2}$ and $K_{2}:=b^{-2}[(1-\delta _{1})k+\delta _{3}-\delta _{2}]$.  In this case, we do not exclude $(\alpha,\beta)$-metrics of Randers type.

\vskip 8mm

\section{Einstein metrics on $S^3$}

We now consider a special family of Randers metrics  $\bar{F}=\alpha+\beta$ of constant flag curvature $K=1$ on $S^3$. This family of Randers metrics were first  introduced in  \cite{BaSh}.
We shall use them to construct Einstein $(\alpha,\beta)$-metrics $F= \alpha \phi(\beta/\alpha)$ with ${\rm Ric}= 2 F$, ${\rm Ric}=0$ and ${\rm Ric}=-2 F$, respectively, but none of them are of constant flag curvature.

We view $S^3$ as a Lie group and let $\eta^1, \eta^2,\eta^3$ be the standard right invariant $1$-form on $S^3$ such that
\[ d \eta^1 =2 \eta^2 \wedge \eta^3, \ \ \ \ d\eta^2 = 2 \eta^3 \wedge \eta^1, \ \ \ \ d\eta^3 = 2\eta^1 \wedge \eta^2.\]
For any number $\e \geq 0$, let $\theta^1: = (1+\e)\eta^1$, $\theta^2:= \sqrt{1+\e}\eta^2$ and $\theta^3:= \sqrt{1+\e}\eta^3$.  Let
\[
\alpha := \sqrt{ [\theta^1]^2 +  [\theta^2]^2 + [\theta^3]^2 }, \ \ \ \ \ \beta := b \theta^1,
\]
where $b=\sqrt{{\e}/{(1+\e)}} <1$. The Levi-Civita connection forms $(\theta_j^{\ i} )$  are given by
\[
d \theta^i = \theta^j \wedge \theta_j^{\ i},
\]
where $\theta_j^{\ i} + \theta_i^{\ j}=0$ and
\[ \theta_2^{\ 1} = \theta^3, \ \ \ \ \theta_3^{\ 1} = -\theta^2, \ \ \ \theta_3^{\ 2} = \frac{1-\epsilon}{1+\epsilon} \theta^1.\]
Then the Riemann curvature tensor of $\alpha$ is given by
\[  d\theta_j^{\ i}-\theta_j^{\ m} \wedge \theta_m^{\ i} = \frac{1}{2} \bar{R}^{\ i}_{j \ pq}\theta^p \wedge \theta^q,\]
where
$\bar{R}^{\ i}_{j \ pq}+ \bar{R}^{\ j}_{i \ pq}=0$, $\bar{R}^{\ i}_{j \ pq}+ \bar{R}^{\ i}_{j \ qp}=0$ and
\begin{eqnarray*}
\frac{1}{2} \bar{R}^{\ 1}_{2 \ pq} \theta^p \wedge \theta^q   & = & \theta^1\wedge \theta^2 ,\\
\frac{1}{2}\bar{R}^{\ 1}_{3 \ pq}  \theta^p \wedge \theta^q & = & \theta^1 \wedge \theta^3 ,\\
\frac{1}{2}\bar{R}^{\ 2}_{1 \ pq} \theta^p \wedge \theta^q  & = & -\theta^1\wedge \theta^2 ,\\
\frac{1}{2}\bar{R}^{\ 2}_{3 \ pq} \theta^p \wedge \theta^q  & = & \lambda  \theta^2 \wedge \theta^3 ,\\
\frac{1}{2}\bar{R}^{\ 3}_{1 \ pq} \theta^p \wedge \theta^q  & = & -\theta^{1}\wedge \theta^3 ,\\
\frac{1}{2}\bar{R}^{\ 3}_{2 \ pq} \theta^p \wedge \theta^q  & = & -\lambda  \theta^2 \wedge \theta^3,
\end{eqnarray*}
where $\lambda := (1-3\e)/(1+\e) = 1-4 b^2$. We get an expression for $\bar{R}^i_{\ j} = \bar{R}^{\ i}_{p \ j q}y^py^q$:
\begin{eqnarray*}
\bar{R}^{1}_{\ 1} & = &  (y^2)^2 +(y^3)^2 ,\\
\bar{R}^{ 2}_{\ 2} & = &  (y^1)^2 + \lambda (y^3)^2 ,\\
\bar{R}^3_{\ 3} & = & (y^1)^2 + \lambda (y^2)^2 ,\\
\bar{R}^1_{\ 2} & = & - y^1y^2 ,\\
\bar{R}^1_{\ 3} & = & -y^1y^3 ,\\
\bar{R}^2_{\ 3} & = & -\lambda y^2y^3 .
\end{eqnarray*}
Then the Ricci curvature of $\alpha$ is given  by
\[ \overline{\rm Ric}= 2 (y^1)^2+(1+\lambda) [ (y^2)^2+(y^3)^2] = 2 \{ \alpha^2 -2( b^2 \alpha^2- \beta^2) \} .\]

By definition, $db_i -b_j \theta_i^{\ j} = b_{i;j} \theta^j.$  A direct computation gives
\[ b_{1;1}=b_{1;2}=b_{1;3}=0,\]
\[ b_{2;1}= b_{2;2}=0, \ \ \ b_{2;3}=-b,\]
\[ b_{3;1}=0,\ \ \ b_{3;2}=b, \ \ \ b_{3;3}=0.\]
We see that
\[
 r_{ij}=0, \ \ \ \  s_j=0.
 \]
We have
\[  s^1_{\ 0} =0, \ \ \ \ s^2_{\ 0} = - b y^3, \ \ \ \ s^3_{\ 0} = b y^2,\]
\[ s_{01}=0, \ \ \ s_{02}= b y^3, \ \ \  s_{03}= - by^2.\]
We have
\[ s^1_{\ m}s^m_{\ 1}=0, \ \ \ s^1_{\ m}s^m_{\ 2}=0, \ \ \ s^1_{\ m}s^m_{\ 3}=0,\]
\[ s^2_{\ m}s^m_{\ 1}=0, \ \ \ s^2_{\ m}s^m_{\ 2}=- b^2, \ \ \ s^2_{\ m}s^m_{\ 3}=0,\]
\[ s^3_{\ m}s^m_{\ 1}=0, \ \ \ s^3_{\ m}s^m_{\ 2}=0, \ \ \ s^3_{\ m}s^m_{\ 3}=- b^2,\]
In short,
\[ s_{0m}s^m_{\ 0} = - (b^2 \alpha^2 - \beta^2 ).\]

By direction computation
\[ s_{12;k}= -b \delta_{2k}, \ \ \ \  s_{13;k} = -b \delta_{3k}, \ \ \ \ s_{23;k}=0.\]
We have
\[  s^1_{\ 0;0}=  - b (y^2)^2-b (y^3)^2,  \ \ \ \ \
s^2_{\ 0;0} =  b y^1y^2, \ \ \ \ s^3_{\ 0;0} =   b y^1 y^3 .\]
\[ s^1_{\ 1;0} =0, \ \ \ \ s^1_{\ 2;0} = - b y^2, \ \ \ \ s^1_{\ 3;0}=-b y^3 \]
\[ s^2_{\ 1;0} =b y^2, \ \ \ \ s^2_{\ 2;0} = 0, \ \ \ \ s^2_{\ 3;0}=0 \]
\[ s^3_{\ 1;0} =b y^3, \ \ \ \ s^3_{\ 2;0} = 0, \ \ \ \ s^3_{\ 3;0}=0 \]
\[ s^1_{\ 0;1}=0, \ \ \ \ s^1_{\ 0;2}= -b y^2, \ \ \ \  s^1_{\ 0;3} = -b y^3,\]
\[ s^2_{\ 0;1} = 0, \ \ \ \ s^2_{\ 0;2} =b y^1, \ \ \ s^2_{\ 0;3} =0,\]
\[ s^3_{\ 0;1} = 0, \ \ \ \  s^3_{\ 0;2}= 0, \ \ \  s^3_{\ 0;3} =b y^1\]
We get
\[ s^m_{\ 0;m} = 2 b y^1 = 2\beta .\]

\bigskip

Therefore  $\alpha$ and $\beta$ satisfy (\ref{EQ2}), (\ref{EQ1}), (\ref{EQ4}) and (\ref{EQ5}) with
\[  \tau = - b^2, \ \ \ \ K_1 =-b^{-2},\ \ \ \  K_2 = 2 b^{-2} .\]

(1) Let $\delta_2 =-b^{-2} (1-\delta_1)$, $\delta_3 = 2$ and $k=-b^{-2}$. Then (\ref{ODE}) is simplified to
\be
(b^2-s^2)[ (1+sQ)Q_s-Q^2 + 2] - [ (2sQ+b^2Q^2)  +1 ] +  \phi^2=0,\label{ODEsphere}
\ee
where $Q:=\phi'/(\phi-s\phi')$.
Note that
$\phi=1+s$ is a special solution of (\ref{ODEsphere}).
For any $\phi =\phi(s)$ satisfying (\ref{ODEsphere}), the metric $F=\alpha \phi(\beta/\alpha)$ is an Einstein metric with ${\rm Ric} = 2 F^2$ on $S^3$.

(2) Let  $\delta_2 =- b^{-2}$ and $\delta_3 =2 -b^{-2} $ and $k=0$. Then (\ref{ODE}) is simplified to
\be
(b^2-s^2)[ (1+sQ)Q_s-Q^2 +2]-[(2sQ+b^2Q^2)+ 1]   = 0. \label{ODEk=0}
\ee
For any $\phi =\phi(s)$ satisfying (\ref{ODEk=0}), the metric $F=\alpha \phi(\beta/\alpha)$ is an Einstein metric  with ${\rm Ric} =0$ on $S^3$.

(3) Let $\delta_2 = - b^{-2}(1-\delta_1)$, $\delta_3 = 2- 2 b^{-2}$ and $k=b^{-2}$. Then (\ref{ODE}) is simplified to
\be
(b^2-s^2)[ (1+sQ)Q_s-Q^2 +2]-[(2sQ+b^2Q^2) +1]  -\phi^2 =0. \label{ODEk=-1}
\ee
For any $\phi =\phi(s)$ satisfying (\ref{ODEk=-1}), the metric $F=\alpha \phi(\beta/\alpha)$ is an Einstein metric with ${\rm Ric} = - 2 F^2$ on $S^3$.

\bigskip
To show that the above Finsler metrics on $S^3$ are not of constant flag curvature, we need to compute $RT^i_{\ j}$ in (\ref{abRiemannCur}).
We  have
\begin{eqnarray*}
 RT^1_{\ 1} & = &  - b \alpha^{-1}[ (y^2)^2+(y^3)^2] \Big \{  C_{311} y^1 - Q_s b\alpha \Big \},\\
RT^1_{\ 2} & = & - b \alpha^{-1} y^2 \Big \{  Q \alpha^2   + [ (y^2)^2+(y^3)^2] C_{311}  \Big \},\\
RT^1_{\ 3} & = &- b \alpha^{-1} y^3 \Big \{  Q \alpha^2   + [ (y^2)^2+(y^3)^2] C_{311}  \Big \} ,\\
RT^2_{\ 1} & = & -b \alpha y^2 [ C_{333} b^2 + Q ]+ b\alpha^{-1}  y^1y^2 \Big \{  C_{311} y^1 - Q_s \alpha b  - b C_{332} \alpha \Big\},  \\\
RT^2_{\ 2} & = & - b^2 [ C_{331} (y^3)^2 + C_{332} (y^2)^2   ]+ b^2 Q^2 \alpha^2 + 2b Q \alpha y^1 + C_{311} b \alpha^{-1} y^1 (y^2)^2 ,\\
RT^2_{\ 3}  & = & b^2 [ C_{331}-C_{332}] y^2y^3 + b \alpha^{-1} C_{311} y^1y^2y^3  ,\\
RT^3_{\ 1} & = & -b^2 C_{332} y^1y^3 - \alpha b^3 C_{333} y^3 - b \alpha Q y^3 + b C_{311}\alpha^{-1} (y^1)^2 y^3 ,\\
RT^3_{\ 2} & = & b^2 [C_{331} - C_{332} ]y^2y^3+ b \alpha^{-1} C_{311} y^1y^2y^3 ,\\
RT^3_{\ 3} & = &   -b^2 [ C_{331} (y^2)^2 + C_{332} (y^3)^2 ] + b^2  Q^2\alpha^2  + 2Q b \alpha y^1 + b \alpha^{-1} C_{311} y^1 (y^3)^2.
\end{eqnarray*}
Further, we obtain
\begin{eqnarray*}
R^1_{\ 1} & = & [  1+ sQ + (b^2-s^2) Q_s ] [ (y^2)^2+(y^3)^2] ,\\
R^2_{\ 2} & = & (y^1)^2+\lambda (y^3)^2 + [b^2Q^2+2sQ]\alpha^2 -[ b^2 C_{332} -  s C_{311}] (y^2)^2-b^2 C_{331} (y^3)^2\\
& = & \frac{ b^2(1+b^2Q^2+2sQ)}{b^2-s^2}[ (y^2)^2+(y^3)^2] -[ b^2 C_{332} -  s C_{311}] (y^2)^2 \\
& & -[b^2 C_{331}-\lambda] (y^3)^2  ,\\
R^3_{\ 3} & = & (y^1)^2+\lambda (y^2)^2 + [b^2Q^2+2sQ]\alpha^2 -[ b^2 C_{332} -  s C_{311}] (y^3)^2-b^2 C_{331} (y^2)^2\\
& = & \frac{ b^2(1+b^2Q^2+2sQ)}{b^2-s^2}[ (y^2)^2+(y^3)^2]  -[ b^2 C_{332} -s C_{311} ](y^3)^2 \\
& & -[  b^2 C_{331}-\lambda] (y^2)^2 .
\end{eqnarray*}
Clearly, $R^1_{\ 1}\not=0$.  Thus when ${\rm Ric}=0$,  we find a Ricci-flat Finsler metric on $S^3$ which is not of zero flag curvature. Also, we can check that $R^{1}_{\ 1}\neq F^{2}(1-l^{1}l_{1})$ and $R^{1}_{\ 1}\neq - F^{2}(1-l^{1}l_{1})$, where $l_{i}:=F_{y^{i}}$ and $l^{i}=g^{ij}l_{j}$.  Thus we have found Einstein Finsler metrics on $S^3$ with $K=1, \ K=0$ and $K=-1$, respectively, which are not of constant flag curvature.

\bigskip

It is surprised that there are singular  Einstein $(\alpha,\beta)$-metrics with ${\rm Ric}=0$ or ${\rm Ric}=-(n-1)F^2$ on $S^3$. This is impossible if the metric is regular by the classical comparison theorems in Finsler geometry.

\section{The conclusion and discussion}

The research of this paper is driven by two motivitions. The first motivition is from S. S. Chern's question, that is, whether or not every smooth manifold admits an Einstein  Finsler metric with Ricci-constant.  We use a Killing form  with non-zero constant length on a Riemannian manifold to construct a family of $(\alpha,\beta)$-metrics and we find equations  that characterize Einstein metrics among this family of $(\alpha,\beta)$-metrics (see Theorem \ref{mainthm}). By Theorem \ref{mainthm}, we know that these constructed  Einstein $(\alpha,\beta)$-metrics may be Einstein Randers metrics or Ricci-flat $(\alpha,\beta)$-metrics of non-Randers type, or  Einstein almost regular $(\alpha,\beta)$-metrics. The second motivition is from the following question: is there any Einstein-Finsler metric which is not of isotropic flag curvature on three dimensional manifolds?  We have found Einstein $(\alpha,\beta)$-metrics on $S^3$ with $Ric =2F, Ric=0$ and $Ric =-2F$, respectively, but none of them are of constant flag curvature (see Theorem \ref{thmn=3}).

Recently, we have learned that L. Huang had constructed a two-parameter family of almost regular Finsler metrics on $S^3$ . His metrics are of constant Ricci curvature $+1$ with  $Ric =2F$, but the flag curvature is  not constant (see \cite{huang}).  However, the technique and method used in \cite{huang} is quite different from ours in this paper.

\end{document}